\documentclass{article}
\usepackage{graphicx}
\usepackage{latexsym}
\usepackage{amsmath}
\usepackage[psamsfonts]{amssymb}

\newtheorem{Theorem}{{\bf Theorem}}
\newtheorem{Lemma}{{\bf Lemma}}
\newtheorem{Corollary}{{\bf Corollary}}
\newtheorem{Definition}{{\bf Definition}}
\newcommand{\RR}{\mathbb R}
\newcommand{\EE}{\mathbb E}

\begin{document}
\title{Statistical Learning Theory of\\ Quasi-Regular Cases}
\author{Koshi Yamada\thanks{Department of Computational Intelligence and Systems Science, Tokyo Institute of Technology, Mail box:G5-19, 4259 Nagatsuta, Midori-ku, Yokohama, 226-8502, Japan\hspace{2mm}\texttt{E-mail:yamada.k.am@m.titech.ac.jp,swatanab@dis.titech.ac.jp}}\and Sumio Watanabe$^{*}$}

\maketitle

\begin{abstract}
Many learning machines such as normal mixtures and layered neural networks
are not regular but singular statistical models, 
because the map from a parameter to a probability distribution is 
not one-to-one. 
The conventional statistical asymptotic theory can not be
applied to such learning machines because the likelihood function
can not be approximated by any normal distribution. 
Recently, new statistical theory has been established 
based on algebraic geometry and it was clarified that
the generalization and training errors are determined by 
two birational invariants, the real log canonical threshold and 
the singular fluctuation. However, their concrete values are
left unknown. 
In the present paper, we propose a new concept, a quasi-regular 
case in statistical learning theory. A quasi-regular case is not
a regular case but a singular case, however, it has the same property as
a regular case. In fact, we prove that, in a quasi-regular case, 
two birational invariants are equal to each other, resulting
that the symmetry of the generalization
and training errors holds. Moreover, the concrete values of two birational
invariants are explicitly obtained, the quasi-regular 
case is useful to study statistical learning theory. 
\end{abstract}

\section{Introduction}

A lot of statistical learning machines which are being applied to
pattern recognition, bioinformatics, robotic control, and artificial 
intelligence have hidden variables, hierarchical layers, and submodules,
because they 
are used to estimate the structure of the true distributions. 
In such learning machines, the map taking parameters
to probability distributions is not one-to-one and the Fisher information
matrices are singular, hence they are called singular learning machines.
For example, three-layered neural networks,
normal mixtures, hidden Markov models, Bayesian networks, and 
reduced rank regressions
are singular learning machines 
\cite{Aoyagi0,Aoyagi,Hagiwara,Hartigan,Hayasaka,1995}. 
If a statistical model is singular, then either 
the maximum likelihood estimator is not subject to the normal distribution even asymptotically
or the Bayes posterior distribution can not be 
approximated by any normal distribution. Hence it has been difficult to 
study their learning performace and to estimate the generalization error from
the training error. 

Recently, new statistical theory has been established based 
on algebraic geometrical method \cite{2001a,Cambridge,2010a,2010b} and it was clarified
that the generalization and training errors in Bayes estimation, 
$G_{n}$ and $T_{n}$, are given by two
birational invariants, the real log canonical threshold $\lambda$ and singular
fluctuation $\nu$ by the formulas, 
\begin{eqnarray}
\EE[G_{n}]&=& \Bigl(\frac{\lambda-\nu}{\beta}+\nu\Bigr)\frac{1}{n}+o(\frac{1}{n}),
\label{eq:Gn}\\
\EE[T_{n}]&=& \Bigl(\frac{\lambda-\nu}{\beta}-\nu\Bigr)\frac{1}{n}+o(\frac{1}{n}),
\label{eq:Tn}
\end{eqnarray}
where $\EE[\;\;]$ shows the expectation value over all training sets, $n$ is the 
number of training samples and $\beta$ is the inverse 
temperature of the Bayes posterior distribution. Based on this relation,
we can define an information criterion which enables us to estimate the
generalization error from the training error \cite{2010a}.

It is well known that, if the true distribution and the statistical model
are in a regular case, then $\lambda=\nu=d/2$ holds where $d$ is the dimension
of the parameter space. In this case, the symmetry of the generalization 
and training errors holds, 
\begin{eqnarray}
\EE[G_{n}]&=& \frac{d}{2n}+o(\frac{1}{n}),
\label{eq:Gn2}\\
\EE[T_{n}]&=& -\frac{d}{2n}+o(\frac{1}{n}),
\label{eq:Tn2}
\end{eqnarray}
for arbitrary $0<\beta\leq \infty$. 
This case corresponds to the well-known Akaike 
Information criterion for regular statistical models. 
However, if they are not in a regular case, neither of 
them is equal to $d/2$ in general. 
Therefore, in order to study singular learning machines, researches 
on two birational invariants are necessary. 

In the present paper, in order to investigate the mathematical 
structure of birational invariants, 
we firstly introduce a new concept, a quasi-regular case,
which satisfies the relation, 
\begin{eqnarray*}
\mbox{Regular} \subsetneq  \mbox{Quasi-Regular} 
\subsetneq  \mbox {Singular}.
\end{eqnarray*}
In other words, a quasi-regular case is not a regular case, however, it has 
the same properties as the regular case. In fact, we prove that,
in quasi-regular cases, both birational invariants are equal to each other,
$\lambda=\nu$, 
and the symmetry of the generalization and training errors holds.
In a quasi-regular case, 
two birational invariants are obtained explicitly, 
hence it is a useful concept in researches of 
statistical learning theory. 

\section{Framework of Bayes Learning}

In this section, we summarize the framework of the 
Bayes learning, and introduce the well-known results.

\subsection{Generalization and Training Errors}

Firstly, we define the generalization and training errors. 
Let $N$, $n$ and $d$ be natural numbers. 
Let $X_{1},X_{2},...,X_{n}$ be
random variables on ${\RR}^{N}$ which are independently subject to 
the same probability density function as $q(x)$. 
Let $p(x|w)$ be a probability density function of $x$ for a 
parameter $w\in W\subset{\RR}^{d}$, where $W$ is a set of parameters. 
The prior distribution is represented by the probability density function
$\varphi(w)$ on $W$. For a given training set 
\[
X^{n}=\{X_{1},X_{2},...,X_{n}\},
\]
the posterior distribution is defined by
\[
p(w|X^{n})=\frac{1}{Z_{n}}\prod_{i=1}^{n}p(X_{i}|w)^{\beta}\varphi(w)dw,
\]
where $0<\beta <\infty$ is the inverse temperature and $Z_n$ is the normalizing constant.
The case $\beta=1$ is most important because it corresponds to the
strict Bayes estimation. The expectation value over the posterior distribution
is denoted by 
\[
\EE_{w}[\;\;]=\int (\;\;)p(w|X^{n})dw.
\]
The predictive distribution is defined by
\[
p(x|X^{n})=\EE_{w}[p(x|w)].
\]
The generalization and training error, $G_{n}$ and $T_{n}$, are
respectively defined by
\begin{eqnarray*}
G_{n}&=& \int q(x)\log\frac{q(x)}{p(x|X^{n})}dx,\\
T_{n}&=&\frac{1}{n}\sum_{i=1}^{n}\log\frac{q(X_{i})}{p(X_{i}|X^{n})}.
\end{eqnarray*}
The generalization error shows the Kullback-Leibler distance from
the true distribution to the estimated distribution. The smaller 
the generalization error is, the better the learning result is.
However, we can not know the generalization error directly, because 
calculation of $G_{n}$ needs
the expectation value over the unknown true distribution $q(x)$. On the
other hand, the training error can be calculated using only training samples,
in practice, as the log likelihood function.
Hence one of the main purposes of statistical learning theory is 
to clarify the mathematical relation between them. 

\subsection{Two Birational Invariants}
Secondly, we define two birational invariants. \\
The Kullback-Leibler distance from the true distribution $q(x)$
to a parametric model $p(x|w)$ is defined by
\[
K(w)=\int q(x)\log\frac{q(x)}{p(x|w)}dx.
\]
Then $K(w)=0$ if and only if $q(x)=p(x|w)$. 
In this paper, we assume that there exists a parameter $w_{0}$ which 
satisfies $q(x)=p(x|w_{0})$ and that 
$K(w)$ is an analytic function of $w$. 

\begin{Definition}{\bf (Real Log Canonical Threshold)} 
The zeta function of statistical learning is defined by
\[
\zeta(z)=\int K(w)^{z}\varphi(w)dw.
\]
Then $\zeta(z)$ is a holomorphic function on the region
$Re(z)>0$, which can be analytically continued to the unique 
meromorphic function on the entire complex plane \cite{2001a}. All poles of 
the zeta function are real, negative, and rational numbers. If its largest 
pole is $(-\lambda)$, then the real log canonical threshold is
defined by $\lambda$. The order of the pole $z=-\lambda$ is 
referred to as a multiplicity $m$.
\end{Definition}

\begin{Definition}{\bf (Singular Fluctuation)}
The functional variance is defined by
\[
V_{n}=\sum_{i=1}^{n}
\{\EE_{w}[(\log p(X_{i}|w))^{2}]
-\EE_{w}[\log p(X_{i}|w)]^{2}\}.
\]
Then it was proved \cite{2010a} that the expectation value
\[
\nu=\frac{\beta}{2}\lim_{n\rightarrow\infty}\EE[V_{n}]
\]
exists. The constant $\nu$ is called the singular fluctuation. 
\end{Definition}

\begin{Theorem}
The expectation values of the generalization and training errors
are given by eq.(\ref{eq:Gn}) and eq.(\ref{eq:Tn}). Therefore
\[
\EE[G_{n}]=\EE[T_{n}]+\frac{2\nu}{n}+o(\frac{1}{n}).
\]
\end{Theorem}
(Proof) This theorem was proved in \cite{Cambridge,2010a}. (Q.E.D.)
\vskip3mm\noindent
{\bf Remarks}. (1) The real log canonical threshold and the
singular fluctuation are invariant under a birational transform
\begin{eqnarray*}
w &= &g(w'), \\
p(x|w)&\mapsto &p(x|g(w')), \\
\varphi(w)&\mapsto & \varphi(g(w'))|g'(w')|,
\end{eqnarray*}
where $|g'(w')|$ is the Jacobian determinant. Such constants are 
called birational invariants. \\
(2) The real log canonical thresholds for
several learning machines were clarified \cite{Aoyagi0,Aoyagi,Yamazaki}
using resolution of singularities. 
However, the singular fluctuation has been left unknown. This paper
provides the first result which clarifies the concrete values of 
singular fluctuation in a singular case. \\
(3) The real log canonical threshold is a well known birational invariant
in algberaic geometry, which plays an important role in higher dimensional
algberaic geometry. The singular fluctuation was found in statistical 
learning theory. 

\subsection{Regular and Singular}

Thirdly, we define regular and singular cases.

\begin{Definition}
A pair of the true distribution and the parametric model, 
$(q(x),p(x|w))$ is called to be in a regular case if and only if
the set $\{w;q(x)=p(x|w)\}$ consists of a single element $w_{0}$ and 
Fisher information matrix
\[
\int \nabla\log p(x|w_{0})(\nabla \log p(x|w_{0}))^{T}q(x)dx
\]
is positive definite. Otherwise, it is called to be in 
a singular case. 
\end{Definition}

For a regular case, the real log canonical threshold and 
the singular fluctuation have been completely clarified. 

\begin{Theorem}
If a pair $(q(x),p(x|w))$ is in a regular case, then
$\lambda=\nu=d/2$, where $d$ is the dimension of the
parameter. 
\end{Theorem}
(Proof) This theorem was proved in \cite{Cambridge}. (Q.E.D.)

\section{Main Results}

In this section, we define a quasi-regular case. This concept 
is firstly proposed by the present paper. Also the main theorem is
introduced. 

\begin{Definition}i{\bf Quasi-Regular Case}j. 
Assume that there exists a parameter $w_{0}\in W^{o}$ such that
$q(x)=p(x|w_{0})$. Without loss of generality, we can assume that
$w_{0}$ is the origin $w_{0}=0$. 
The original parameter is denoted by $w=(w_{1},w_{2},...,w_{d})$. 
Let $g$ and $\Delta d_{1},\Delta d_{2},...,\Delta d_{g}$ be natural numbers which satisfy
\[
\Delta d_{1}+\Delta d_{2}+\cdots+\Delta d_{g}=d
\]
and $\Delta  d_{0}=0$. 
We define
\[d_j=\Delta d_{0}+\cdots+\Delta d_{j}\: \: (j=0,\cdots,g)\]
and a function  $u=(u_{1},u_{2},...,u_{g})\in{\RR}^{g}$ of
the paramater $w\in{\RR}^{d}$ by 
\begin{eqnarray*}
\displaystyle u_{1}&=& \prod_{j=1}^{d_{1}}w_{j},\\
\displaystyle u_{2}&=& \prod_{j=d_{1}+1}^{d_{2}}w_{j},\\
\cdots &=& \cdots, \\
\displaystyle u_{g}&=& \prod_{j=d_{g-1}+1}^{d_{g}}w_{j}.
\end{eqnarray*}
If there exist constants $c_{1},c_{2}>0$ such that,
for arbitrary $w\in W$, 
\[
c_{1}(u_{1}^{2}+\cdots+u_{g}^{2})
\leq K(w)\leq 
c_{2}(u_{1}^{2}+\cdots+u_{g}^{2}),
\]
then the pair $(q(x),p(x|w))$ is called 
to be in a quasi-regular case.
\end{Definition}

\vskip3mm\noindent
{\bf Remark.} (1) If $g=d$, then 
\[
\{w;q(x)=p(x|w)\}=\{0\}
\]
and the quasi-regular case
corresponds to the regular case. Hence a quasi-regular case 
contains a regular case as a special one. \\
(2) If $d\neq g$, then ``$K(w)=0\Longleftrightarrow w=0$" does not
hold, because, for at least one variable $w_{j}$, $K(0,0,..,w_{j},0,..,0)=0$.
Hence a quasi-regular case with $d\neq g$ is not a regular case but a
singular case. \\
(3) There are singular cases which are not contained in quasi-regular cases. Therefore,
\begin{eqnarray*}
\mbox{Regular} \subsetneq  \mbox{Quasi-Regular} 
\subsetneq  \mbox {Singular}
\end{eqnarray*}
holds. The present paper shows in Theorem \ref{Theorem:Gn} that a quasi-regular case is not 
a regular case, however, it has the same property as a regular case. 
\vskip3mm\noindent
{\bf Example.1} Let a statistical model be 
\[
p(x,y|w)=\frac{r(x)}{\sqrt{2\pi}}
\exp(-\frac{1}{2}(y-ax^{2}-b\tanh(cx))^{2}),
\]
where $w=(a,b,c)$ is the parameter and $r(x)$ is the probability density function  of $x$.
If the true distribution is 
given by $q(x,y)=p(x,y|0,0,0)$, then by using 
\[
u_{1}=a,\;\;u_{2}=bc,
\]
it follows that 
\[
K(w)=
\frac{1}{2}\int (ax^{2}+b\tanh(cx))^{2}r(x)dx
\]
satisfies the condition for a quasi-regular case with $g=2$, 
because $x^{2}$ and $\tanh(cx)/c$ is linearly independent. 
In fact there exist $c_{1},c_{2}>0$ such that 
\[
c_{1}(a^{2}+(bc)^{2})\leq K(w)
\leq c_{2}(a^{2}+(bc)^{2}).
\]
Hence the set of true parameters consists of 
the union of two lines, 
\[
\{w;q(x,y)=p(x,y|w)\}= \{a=0, bc=0\}.
\]
\vskip3mm\noindent
{\bf Example.2} Let a statistical model be 
\[
p(x,y|w)=\frac{r(x)}{\sqrt{2\pi}}
\exp(-\frac{1}{2}(y-ax-b\tanh(cx))^{2}),
\]
where $w=(a,b,c)$ is the parameter and the true distribution is 
given by $q(x,y)=p(x,y|0,0,0)$. Then 
because $x$ and $\tanh(cx)/c$ is not linearly independent
as $c\rightarrow 0$, hence this case does not satisfies the
quasi-regular condition. In this case
\[
c_{1}((a+bc)^{2}+b^{2}c^{6})\leq K(w)
\leq c_{2}((a+bc)^{2}+b^{2}c^{6}). 
\]
Example.2 resembles Example.1, however, from the
viewpoint of statistical learning theory, they are
different. 
\vskip3mm\noindent
{\bf Example.3} 
Let a statistical model be 
\begin{eqnarray*}
p(x,y,z|w)=\frac{r(x,y)}{\sqrt{2\pi}}
\exp(-\frac{1}{2}(z-f(x,y,w))^{2}),
\end{eqnarray*}
where 
\begin{eqnarray*}
f(x,y,w)&=&
a_{1}\sin(b_{1}x)+a_{2}x\sin(b_{2}x)\\
&& + a_{3}\sin(b_{3}y)+a_{4}y\sin(b_{4}y),
\end{eqnarray*}
and 
$w=\{(a_{i},b_{i})\}$ is the parameter and the true distribution is 
given by $q(x,y,z)=p(x,y,z|0)$. Then $(q(x,y,z),p(x,y,z|w))$ is in a 
quasi-regular case with $g=4$. 
\vskip3mm
The following is the main theorem of the present paper. 

\begin{Theorem}{\bf (Main Theorem)}.\label{Theorem:Gn}
Assume that the pair $(q(x),p(x|w))$ is in a quasi-regular
case and that $\varphi(w)>0$ on $W$. Then the real log 
canonical threshold and the singular fluctuation are given by 
\begin{eqnarray*}
\lambda=\nu=\frac{g}{2}
\end{eqnarray*}
and
\begin{eqnarray*}
m = d-g+1.
\end{eqnarray*}
\end{Theorem}

\begin{Corollary}
Assume that the pair $(q(x),p(x|w))$ is in a quasi-regular
case and that $\varphi(w)>0$ on $W$. 
For arbitrary $0<\beta<\infty$ the symmetry of the
generalization and training errors holds, 
\begin{eqnarray*}
\EE[G_{n}]&=& \frac{g}{2n}+o(\frac{1}{n}), \\
\EE[T_{n}]&=& - \frac{g}{2n}+o(\frac{1}{n}). 
\end{eqnarray*}
\end{Corollary}

\vskip3mm\noindent
{\bf Remarks.}(1) The above theorem shows the generalization and training 
errors for Bayes estimation. In the quasi-regular case, 
they have the same property as those in 
regular cases, however, the generalization and training errors of
the maximum likelihood estimation is different from regular case in general. \\
(2) In the maximum likelihood method, 
the training error of a singular case is far smaller than that 
of a regular case, whereas the generalization error of a singular case is far 
larger than that of a regular case. From the viewpoint of the maximum likelihood
method, the quasi-regualr case is contained in the singular case. 
In the present paper, we prove that the quasi-regular case has the same 
property as the regular case from the viewpoint of the Bayes estimation.

\section{Proofs}

In this section, we prove the main theorem. At first, we derive the real
log canonical threshold of the quasi-regular case. 

\begin{Lemma}
The real log canonical threshold and its order are given by
$\lambda=g/2$ and $m=d-g+1$ respectively. 
\end{Lemma}
(Proof) Since each function $\{u_{j};j=1,2,...,g\}$ 
does not have common variable $w_{k}$, the real log canonical threshold
is given by the sum of individual real log canonical thresholds 
(Remark 7.2 in \cite{Cambridge}) defined by 
\begin{eqnarray*}
\zeta_{j}(z)&=& \int \prod_{i=d_{j-1}+1}^{d_{j}}(w_{i})^{2z}dw_{i}\\
&=& \frac{C}{(z+1/2)^{d_{j}-d_{j-1}}}+\cdots+.
\end{eqnarray*}
Hence $\lambda$ is equal to $g$ times $1/2$, hence 
$\lambda=g/2$. The multiplicity is also given by 
\begin{eqnarray*}
m&=&d_{1}+d_{2}-d_{1}+\cdots +d_{g}-d_{g-1}-(g-1)\\
&=& d-g+1,
\end{eqnarray*}
which shows the Lemma. (Q.E.D.) 

\begin{Definition}
For a given pair of the true distribution $q(x)$ and 
the parametric model $p(x|w)$, the log density
ratio function is defined by
\[
f(x,w)=\log\frac{q(x)}{p(x|w)}.
\]
\end{Definition}

The following lemma shows that the log density ratio
function of the quasi-regular case is represented by 
$g$ linearly independent functions. 

\begin{Lemma} \label{Lemma:ideal}
Assume that the pair $(q(x),p(x|w))$ is 
in a quasi-regular case. Then there exists a set of
functions $\{e_{j}(x,u);j=1,2,...,g\}$ which are 
analytic functions of $u$ and 
\[
f(x,w)=\sum_{j=1}^{g}u_{j}e_{j}(x,u)
\]
in an open neighborhood of $u=0$. 
\end{Lemma}
(Proof) Let us define a function
\[
F(t)=t+e^{-t}-1.
\]
for $t\in {\RR}^{1}$. Then $F(0)=0$, $F'(0)=0$, and
$F''(0)=1$, resulting that $F(t)\geq 0$ and that
$F(t)=0$ if and only if $t=0$. 
Moreover, $F(t)\cong (1/2)t^{2}$ for small $|t|$. 
Therefore, 
\begin{eqnarray}
K(w)&=&\int q(x)F\Bigl(\log\frac{q(x)}{p(x|w)}\Bigr)dx \nonumber \\
&=& \int q(x) F(f(x,w))dx\nonumber \\
&\cong & \frac{1}{2}\int q(x) f(x,w)^{2}dx.\label{eq:f2}
\end{eqnarray}
By the assumption of the quasi-regular case, 
$K(w)=0$ if and only if $u_{1}=u_{2}=\cdots = u_{g}=0$,
which is equivalent to $f(x,w)\equiv 0$. 
That is to say, $f(x,w)$ is contained in the ideal 
of analytic functions 
generated by $u_{1},u_{2},...,u_{g}$. Hence there exist 
a set  $\{e_{j}(x,u)\}$ of analytic functions of $u$, which 
satisfies 
\[
f(x,w)=\sum_{j=1}^{g}u_{j}e_{j}(x,u).
\]
Therefore, we obtained the Lemma. (Q.E.D.)
\vskip3mm\noindent
In the following lemma, we show that the quasi-regular case
has the generalized Fisher information matrix. 

\begin{Lemma}
The $g\times g$ matrix $I(u)$ is
defined by
\[
I_{ij}(u)\equiv \int q(x)e_{i}(x,u)e_{j}(x,u)dx.
\]
Then $I(u)$ 
is positive definite in an open neighborhood of $u=0$. 
\end{Lemma}
(Proof)  By Lemma
\ref{Lemma:ideal} and eq.(\ref{eq:f2}),
in the neighborhood of $u=0$,
\[
K(w)=\frac{1}{2}(u \cdot I(u) u).
\]
By the condition of the quasi-regular case, 
\[
c_{1}\sum_{j=1}^{g}u_{j}^{2}\leq K(w).
\]
Hence the minimum eigenvalue of $I(u)$ is positive, which shows 
$I(u)$ is positive definite. (Q.E.D.)
\vskip3mm\noindent
The following definition and lemma show that 
the empirical loss function of the quasi-regular case 
has the same decomposition as that of the regular case. 

\begin{Definition}
A random process $\xi_{n}(u)\in{\RR}^{g}$ is defined by
\[
\xi_{n}(u)=\frac{1}{\sqrt{n}}\sum_{i=1}^{n}
\{\frac{1}{2}I(u)u-e(X_{i},u)\}
\]
where
\[
e(x,u)=(e_{1}(x,u),e_{2}(x,u),...,e_{g}(x,u))^{T}.
\]
\end{Definition}

\begin{Lemma}
The empirical loss function defined by
\[
K_{n}(w)=\frac{1}{n}\sum_{i=1}^{n}f(X_{i},w)
\]
is represented by 
\[
K_{n}(w)=\frac{1}{2}(u,I(u)u)-\frac{1}{\sqrt{n}}\;u\cdot \xi_{n}(u)
\]
in the neighborhood of $u=0$. Moreover, the random process $\xi_{n}(u)$ 
converges to the gaussian process $\xi(u)$ that satisfies
\[
\EE[\xi(0)\cdot I(0)^{-1}\xi(0)]=g.
\]
\end{Lemma}
(Proof) The empirical loss function is given by
\[
K_{n}(w)=K(w)-\frac{1}{n}
\sum_{i=1}^{n}\{K(w)-f(X_{i},w)\}.
\]
By combining this equation with 
the definition of $\xi_{n}(u)$, the first half of the Lemma is
obtained.  For the second half, the convergence $\xi_{n}(u)$ is 
derived from the general empirical process theory. Moreover, 
\begin{eqnarray*}
&& \EE[\xi_{n}(0)\cdot I(0)^{-1}\xi_{n}(0)]\\
&&= \EE[\mathrm{tr}(I(0)^{-1}\xi_{n}(0)\xi_{n}(0)^{T})]=g,
\end{eqnarray*}
where we used the covariance matrix of $\xi_{n}(0)$ 
\[
\EE[\xi_{n}(0)\xi_{n}(0)^{T}]=\int q(x)e(x,0)e(x,0)^{T}dx
=I(0),
\]
which completes the Lemma. 
(Q.E.D.)
\vskip3mm\noindent
In the quasi-regular case, the relation between 
$w=(w_{1},w_{2},...,w_{d})$ and $u=(u_{1},u_{2},...,u_{g})$
is important. The following lemma shows the property of
the quasi-regular case. This lemma does not hold in general
singular cases. 

\begin{Lemma}
When $n$ tends to infinity, 
\begin{eqnarray*}
\prod_{j=1}^{g}\delta \Bigl(\frac{u_{j}}{\sqrt{n}}
-\prod_{k=d_{j-1}+1}^{d_{j}}w_{k}\Bigr)
\cong c_{3}(\log n)^{m-1}\prod_{j=1}^{d}
\delta(w_{j})
\end{eqnarray*}
where $m=d-g+1$ and $c_{3}>0$ is a constant. 
\end{Lemma}
(Proof) Firstly, we prove that 
the delta function with variables ${\bf x}=(x_{1},x_{2},...,x_{d})$ 
in $M\equiv[0,1]^{d}$
\[
D(t,{\bf x})=\delta(t-x_{1}x_{2}\cdots x_{d})
\]
has asymptotic expansion for $t\rightarrow 0$, 
\begin{eqnarray}
D(t,{\bf x})
&=& \frac{ (-\log t)^{d-1}}{(d-1)!} \prod_{k=1}^{d}\delta(x_{k})
\nonumber \\
&&+o((-\log t)^{d-2}).\label{eq:delta}
\end{eqnarray}
Let $\phi({\bf x})$ be an arbitrary $C^{\infty}$-class function of ${\bf x}$
whose support is containd in 
\[
D_{t}(\phi)\equiv \int_{M} D(t,{\bf x})\phi({\bf x})d{\bf x}.
\]
Then its Mellin transform is
\[
\int D_{t}(\phi)t^{z}dt =\int_{M} \prod_{i=1}^{d}(x_{i})^{z}\phi({\bf x})d{\bf x},
\]
where $M$ is the compact set that is the support of $\phi$. 
Without loss of generality 
By using Taylor expansion
\[
\phi({\bf x})=\phi(0)+{\bf x}\cdot \nabla\phi(0)+\cdots,
\]
we have the asymptotic expansion,
\[
\int D_{t}(\phi)t^{z}dt=\frac{1}{(z+1)^{d}}\phi(0)+\cdots.
\]
Therefore 
\[
\int D_{t}(t,{\bf x})t^{z}dt =\frac{1}{(z+1)^{d}}\prod_{k=1}^{d}\delta(x_{k})
+\cdots
\]
for ${\bf x}\in [0,1]^{d}$. By using inverse Mellin transform, we
obtained eq.(\ref{eq:delta}). 
Secondly, let us prove the Lemma. 
By using eq.(\ref{eq:delta}), for each $u_{j}$, 
\begin{eqnarray*}
&& \delta (\frac{u_{j}}{\sqrt{n}}
-\prod_{k=d_{j-1}+1}^{d_{j}}w_{k})\\
&& \propto 
(\log n)^{d_{j}-d_{j-1}-1}\prod_{j=d_{j-1}+1}^{d_{j}}
\delta(w_{j})
\end{eqnarray*}
when $n\rightarrow\infty$. 
By summing up these relations 
for $j=1,2,...,g$, Lemma is obtained. (Q.E.D.)
\vskip3mm\noindent
Let us return to the proof of the Main theorem. 
\vskip3mm
\noindent
{\bf (Proof of Main Theorem)} It was proved by eq.(6.4) in \cite{Cambridge}
that the expectation value of $K_{n}(w)$ is given by two birational invariants, 
\[
\EE[\EE_{w}[K_{n}(w)]]=\frac{\lambda}{n\beta}-\frac{\nu}{n}+o(\frac{1}{n}).
\]
Since we have already obtained the value of $\lambda$ in Lemma1,
that is to say, $\lambda=g/2$, we can derive 
the value of $\nu$ by calculating $\EE[\EE_{w}[K_{n}(w)]]$.
The posterior distribution is represented by the empirical loss function by
\[
p(w|X^{n})\propto\exp(-n\beta K_{n}(w))\varphi(w)dw.
\]
The integration of the outside of the neighborhood of $u=0$ with respect to
the posterior distribution goes to zero 
with the smaller order than $\exp(-\sqrt{n})$ as Lemma 6.3 in \cite{Cambridge}, 
hence we can restrict the integrated region to the neighborhood of $u=0$. 
The empirical loss function is rewritten as 
\begin{eqnarray*}
K_{n}(w)&=& \frac{1}{2}\|I(u)^{\frac{1}{2}}\Bigl(u-I(u)^{-1}\frac{\xi_{n}(u)}{\sqrt{n}}\Bigr)\|^{2}\\
& & -\frac{1}{2n}(\xi_{n}(u)\cdot I(u)^{-1}\xi_{n}(u)).
\end{eqnarray*}
In the neighborhood of $u=0$, we obtain 
\begin{eqnarray*}
K_{n}(w)&\cong& \frac{1}{2}\|I(0)^{\frac{1}{2}}
\Bigl(u-I(0)^{-1}\frac{\xi_{n}(0)}{\sqrt{n}}\Bigr)\|^{2}\\
& & -\frac{1}{2n}(\xi_{n}(0)\cdot I(0)^{-1}\xi_{n}(0)).
\end{eqnarray*}
For an arbitrary function $F(\;\;)$, 
\begin{eqnarray*}
&& \int F(\sqrt{n}\;u) dw \\
&& = \int F(\sqrt{n}\;u) \prod_{j=1}^{g}\delta \Bigl(u
-\prod_{k=d_{j-1}+1}^{d_{j}}w_{k}\Bigr)dw du\\
&& =\int F(u) \prod_{j=1}^{g}\delta \Bigl(\frac{u}{\sqrt{n}}
-\prod_{k=d_{j-1}+1}^{d_{j}}w_{k}\Bigr)dw \frac{du}{n^{g/2}}\\
&& =\frac{c_{3}(\log n)^{m-1}}{n^{g/2}}\int F(u) du.
\end{eqnarray*}
On the other hand,
\begin{eqnarray*}
nK_{n}(w)&=& \frac{1}{2}\|I(0)^{1/2}
\Bigl(\sqrt{n}\;u-I(0)^{-1}\xi_{n}(0)\Bigr)\|^{2}\\
& & -\frac{1}{2}(\xi_{n}(0)\cdot I(0)^{-1}\xi_{n}(0))\\
&\equiv & \hat{K}_{n}(\sqrt{n}\;u).
\end{eqnarray*}
Therefore, 
\begin{eqnarray*}
\begin{split}
&\EE_{w}[K_{n}(w)]=\frac{\int K_{n}(w)\exp(-n\beta K_{n}(w))\varphi(w)dw}{\int \exp(-n\beta K_{n}(w))\varphi(w)dw}\\
&=\frac{1}{n}\frac{\int \hat{K}_{n}(\sqrt{n}\;u)\exp(-\beta\hat{K}_{n}(\sqrt{n}\;u))\varphi(w)dw}{\int \exp(-\beta\hat{K}_{n}(\sqrt{n}\;u))\varphi(w)dw}\\
&=\frac{1}{n}\frac{\int \hat{K}_{n}(u)\exp(-\beta\hat{K}_{n}(u))du}{\int \exp(-\beta\hat{K}_{n}(u))du}\\
&= \frac{1}{2n}\frac{\int \|I(0)^{\frac{1}{2}}(u-\xi_{n}^{*})\|^{2}\exp(-\beta\hat{K}_{n}(u))du}{\int \exp(-\beta\hat{K}_{n}(u))du}{}\nonumber\\
&{}-\frac{1}{2n}(\xi_{n}(0)\cdot I(0)^{-1}\xi_{n}(0)),
\end{split}
\end{eqnarray*}
where the notation
\[
\xi_{n}^{*}=I(0)^{-1}\xi_{n}(0)
\]
is used. Finally, by the integral formlula
\[
\frac{ 
\int \|I(0)^{1/2}u\|^{2}
\exp(-\frac{\beta}{2}\|I(0)^{1/2}u\|^{2})du
}{
\int
\exp(-\frac{\beta}{2}\|I(0)^{1/2}u\|^{2})du
}
= \frac{g}{\beta}
\]
and by Lemma.4, we have 
\begin{eqnarray*}
\EE[\EE_{w}[K_{n}(w)]] = \frac{g}{2\beta n}-\frac{g}{2n}+o(\frac{1}{n}),
\end{eqnarray*}
Then, because $\lambda=\frac{g}{2}$ holds from Lemma1,  we obtain the Theorem. 
(Q.E.D.)
\vskip3mm\noindent
{\bf Example.4} By the main theorem of this paper, 
the real log canonical threshold and the singular fluctuation
of Example.1 are $\lambda=\nu=1$. 
Also those of Example.3 are $\lambda=\nu=2$. 

\section{Discusion}

Let us discuss the result of this paper from the 
two different points of view. Firstly, we study the 
theoretical aspect and then the practical aspect. 

\subsection{Theoretical point of view}

In the present paper, we introduced a new concept, a quasi-regular case.
A quasi-regular case is not a regular case, but it has the same property as 
the regular case. Table.1 shows comparison of the real log canonical threshold (RLCT),
singular fluctuation (SF), the generalization error $G_{n}$, and the training error $T_{n}$. 

Even for the general singular cases, real log canonical thresholds have been clarified in several cases. 
However, this paper is the first case in which the singular fluctuation was clarified. 
In general singular cases, it is conjectured that the real log canonical threshold is not 
equal to the singular fluctuation. To clarify such conjecture is the future study.

\begin{table}[tb]
\begin{center}
\begin{tabular}{|c|c|c|c|}
\hline
           & regular & quasi-regular &singular \\
\hline
RCLT & $d/2$ & $g/2$ & $ \lambda $ \\
\hline
SF & $d/2$ &  $g/2$ & $\nu$ \\
\hline
$G_{n}$ & $d/(2n)$  & $g/(2n)$ & $ ((\lambda-\nu)/\beta +\nu)/n $ \\
\hline
$T_{n}$ & $-d/(2n)$ & $-g/(2n)$ & $ ((\lambda-\nu)/\beta - \nu)/n $ \\
\hline
\end{tabular}
\end{center}
\caption{Regular, Quasi-Regular, and Singular}
\end{table}

\subsection{Practical point of view}

In applications, even if both birational invariants are unknown, 
the generalization error can be estimated 
from the training error and the functional variance \cite{2010a} because
\[
\EE[G_{n}]=\EE[T_{n}]+\frac{\beta}{n}\EE[V_{n}]+o(\frac{1}{n}),
\]
which is asymptotically equivalent to Bayes cross validation \cite{2010b}. 

However, in Bayes estimation, the method how to 
approximate the posterior distribution using Markov 
chain Monte Carlo (MCMC) method is an important issue.
There are a lot of parameters which determine the MCMC process, 
for example, times of burn-in, times of sufficiently updates, and so on. 
If we know the concrete values of birational invariants,
then we can evaluate how accurate the MCMC process is \cite{Nagata}.
Therefore, the quasi-regular cases are appropriate for evaluating 
MCMC process. 
It is the future study to evaluate MCMC process 
using the quasi-regular cases.

\section{Conclusion}

In the present paper, 
a new concept, a quasi-regular case, was firstly proposed,
and its theoretical foundation was constructed. 
A quasi-regular case is not a regular case but a singular case, 
whereas it has the same property as a regular case. 
In a quasi-regular case, it was proved that
the real log canonical threshold is equal to the singular 
fluctuation. This is the first case in which nontrivial
value of singular fluctuation is clarified. 

\subsection*{Acknowledgement}
This research was partially supported by the Ministry of Education,
Science, Sports and Culture in Japan, Grant-in-Aid for Scientific
Research 23500172.


\begin{thebibliography}{99}


\bibitem{Aoyagi0}
M.Aoyagi,S.Watanabe,``Resolution of singularities and generalization error with 
Bayesian estimation for layered neural network," 
Vol.J88-D-II, No.10, pp.2112-2124, 2005. 

\bibitem{Aoyagi}
M.Aoyagi, S.Watanabe,``Stochastic complexities of reduced rank regression in Bayesian estimation," Neural Networks, Vol.18,No.7, pp.924-933, 2005. 

\bibitem{Atiyah}
M.F.Atiyah,``Resolution of singularities and division of distributions,"
Comm. Pure Appl. Math., Vol.13,pp.145-150,1970. 

\bibitem{Hagiwara} K. Hagiwara,
``On the Problem in Model Selection of Neural Network Regression 
in Overrealizable Scenario,"
Neural Comput., Vol.14,Vol.8, pp.1979 - 2002, 2002.

\bibitem{Hartigan}
J.A.Hartigan,``A failure of likelihood asymptotics for normal mixture,"
 Proc. of Barkeley Conf. in honor of Jerzy Neyman and Jack Keifer, Vol.2,
 pp.807-810,1985.

\bibitem{Hayasaka}
T. Hayasaka, M. Kitahara, and S. Usui,
``On the Asymptotic Distribution of the Least-Squares Estimators 
in Unidentifiable Models," Neural Comput., Vol.16 ,No.1, pp.99 - 114, 2004.

\bibitem{Hironaka}
H. Hironaka, ``Resolution of singularities of an algebraic variety over a field of characteristic zero," Ann. of Math., Vol.79, 109-326,1964.

\bibitem{Kashiwara}
M. Kashiwara, ``B-functions and holonomic systems," 
Inventions Math., 38, 33-53.1976.

\bibitem{Nagata}
K. Nagata, S. Watanabe,
``Asymptotic Behavior of Exchange Ratio in Exchange Monte Carlo Method,''
International Journal of Neural Networks, Vol. 21, No. 7, pp. 980-988, 2008.

\bibitem{1995}
S. Watanabe, ``Generalized Bayesian framework for neural networks with singular Fisher information matrices," Proc. of International Symposium on Nonlinear Theory and Its applications, (Las Vegas), pp.207-210, 1995. 

\bibitem{2001a}
S. Watanabe, "Algebraic Analysis for Nonidentifiable Learning Machines," Neural Computation, Vol.13, No.4, pp.899-933, 2001.

\bibitem{2001b}
S. Watanabe, "Algebraic geometrical methods for hierarchical learning machines," Neural Networks, Vol.14, No.8,pp.1049-1060, 2001. 

\bibitem{2010a}
Sumio Watanabe, "Equations of states in singular statistical estimation", 
Neural Networks, Vol.23, No.1, pp.20-34, 2010. 

\bibitem{2010b}
Sumio Watanabe, ``Asymptotic Equivalence of Bayes Cross Validation and 
Widely Applicable Information Criterion in Singular Learning Theory," 
Journal of Machine Learning Research, Vol.11, (DEC), pp.3571-3591, 2010. 

\bibitem{Cambridge}
S. Watanabe, ``Algebraic geometry and statistical learning theory,"
Cambirdge University Press, 2009. 

\bibitem{Yamazaki}
K.Yamazaki, S.Watanabe,``Singularities in mixture models and upper bounds of stochastic complexity." International Journal of Neural Networks, Vol.16, No.7, pp.1029-1038,2003. 

\end{thebibliography}
\end{document}